\documentclass[11pt]{amsart}

\usepackage{amssymb}

\theoremstyle{plain}
\newtheorem{theorem}{Theorem}[section]
\newtheorem{proposition}{Proposition}[section]

\theoremstyle{definition}
\newtheorem{definition}{Definition}[section]

\newtheorem{remark}{Remark}[section]

\newcommand{\sn}{{\mathbb S}}
\newcommand{\tn}{{\mathbb T}}
\newcommand{\rn}{{\mathbb R}}

\newcommand{\h}{{\mathrm H}}
\renewcommand{\v}{{\mathrm V}}

\DeclareMathOperator{\trace}{trace}
\DeclareMathOperator{\grad}{grad}
\DeclareMathOperator{\ricci}{Ricci}
\DeclareMathOperator{\riem}{Riem}
\DeclareMathOperator{\Div}{div}

\begin{document}
\title{Harmonic maps and Kaluza-Klein metrics on spheres}

\author{M. Benyounes}
\address{D{\'e}partement de Math{\'e}matiques \\
Universit{\'e} de Bretagne Occidentale \\
6, avenue Victor Le Gorgeu \\
CS 93837, 29238 Brest Cedex 3, France}
\email{Michele.Benyounes@univ-brest.fr {\rm{and}} Eric.Loubeau@univ-brest.fr}

\author{E. Loubeau}

\author{L.~Todjihounde}
\address{Institut de Math{\'e}matiques et de Sciences Physiques\\
B.P. 613 Porto-Novo, R{\'e}publique du B{\'e}nin}
\email{leonardt@imsp-uac.org}



\begin{abstract}
This article studies the harmonicity of vector fields on Riemannian manifolds, viewed as maps into the tangent bundle equipped with a family of Riemannian metrics. Geometric and topological rigidity conditions are obtained, especially for surfaces and vector fields of constant norm, and existence is proved on two-tori. Classifications are given for conformal, quadratic and Killing vector fields on spheres. Finally, the class of metric considered on the tangent bundle is enlarged, permitting new vector fields to become harmonic.
\end{abstract}

\maketitle

\section{Introduction}

Though very interesting in many settings, the theory of harmonic maps fails to produce any worthwhile result when applied to vector fields, 
seen as maps from a Riemannian manifold $(M,g)$ into its tangent bundle $TM$ equipped with its simplest metric, the Sasaki metric. This situation has led to consider constrained problems on the same functional, e.g. harmonic sections and harmonic unit sections. However, recently, new classes of metric on $TM$ have been shown to allow a richer existence theory and, with respect to adequate metrics, standard vector fields can produce new harmonic maps, for example a two-parameter family including the Sasaki metric (\cite{BLW1}), $g$-natural metrics (\cite{ACP}) or an ad-hoc Riemannian metric based on a deformation of the horizontal distribution (\cite{Oniciuc}).

The main difficulty here is to strike a balance between the harmonicity of the vector fields and the geometric relevance of the metric on $TM$. In this paper, given a Riemannian manifold $(M,g)$, we consider on $TM$ Riemannian metrics in the intersection of the largest known class of metrics on tangent bundles, i.e. $g$-natural metrics, and Kaluza-Klein metrics, as commonly defined on principal bundles (cf.~\cite{Wood}). 

Recall that, at a point  $(p,e) \in TM$, the tangent space $T_{(p,e)}TM$ splits into its horizontal and vertical spaces (\cite{Dombrowski}):
\begin{equation}\label{eq1}
 T_{(p,e)}TM = \h_{(p,e)} \oplus \v_{(p,e)} ,
\end{equation}
where $\v_{(p,e)}$ is the kernel of the differential of the canonical projection $\pi : TM \to M$ and $\h_{(p,e)}$ is the kernel of the connecting map 
\[
 K_{(p,e)} =K : T_{(p,e)}TM \to T_p M ,  \quad 
 K(V) = d (\exp_p \circ R_{-e} \circ\tau)(V),
\]
where $\tau : U \subset TM \to T_p M$ sends a vector 
$v\in T_q M$, with $(q,v) \in U$, $U$ being an open neighbourhood 
of $(p,e)$ in $TM$, by parallel transport along the unique
geodesic from $q$ to $p$, to a vector in $T_p M$ and the map $R_{-e}$ is simply the translation by $-e$ in $T_p M$. 
One can check that $\h_{(p,e)}\cap \v_{(p,e)} = \{0\}$ 
and $\h_{(p,e)} \oplus \v_{(p,e)}=T_{(p,e)}TM$.
Any vector in $T_{(p,e)}TM$ can be decomposed into its horizontal and vertical parts and any vector $X \in T_p M$ admits a horizontal lift $X^h \in\h_{(p, e)}$ and a vertical lift $X^v \in\v_{(p, e)}$ defined by
\[ 
 K_{(p,e)}(X^v) = X,  \quad  
d\pi_{(p,e)}(X^h ) = X.
\]
Metrics on $TM$ can therefore be characterized by their values on horizontal and vertical lifts.
\begin{definition}~\cite{AS}
Let $(M,g)$ be a Riemannian manifold, a metric $G$ on $TM$ will be called $g$-natural if, at the point $(p,e)\in TM$, it has the form
\begin{align*}
G(X^h , Y^h) &= A(|e|^2)g(X,Y) + D(|e|^2)g(X,e)g(e,Y) ;\\
G(X^h , Y^v) &= E(|e|^2)g(X,Y) + F(|e|^2)g(X,e)g(e,Y) ;\\
G(X^v , Y^v) &= B(|e|^2)g(X,Y) + C(|e|^2)g(X,e)g(e,Y) ,
\end{align*}
where $A,B,C,D,E$ and $F$ are real $C^2$-functions of $|e|^2$. Conditions are needed on $A,B,C$ and $D$ to ensure that $G$ is positive definite.
\end{definition}

\begin{remark}
To have $\pi : TM \to M$ conformal submersion, we need $D= 0$ and to have the horizontal and the vertical distributions orthogonal one to the other, we need $E= F= 0$.
\end{remark}

The class of metrics we study in this article sits inside $g$-natural metrics but retains the geometric properties of the tangent bundle.

\begin{definition}
Let $(M,g)$ be a Riemannian manifold, a metric $G$ on $TM$ will be called Kaluza-Klein if, at the point $(p,e)\in TM$, it takes the form
\begin{align*}
G(X^h , Y^h) &= A(|e|^2)g(X,Y)  ;\\
G(X^h , Y^v) &= 0 ;\\
G(X^v , Y^v) &= B(|e|^2)g(X,Y) + C(|e|^2)g(X,e)g(e,Y) ,
\end{align*}
where $A,B$ and $C$ are real functions of $|e|^2$. The functions $A$ will be assumed strictly positive and $B$ and $C$ such that $G$ is positive definite.
\end{definition}

Standard computations with the Koszul formula or general results for $g$-natural metrics from~\cite{AS}, give the expression of the Levi-Civita connection.

\begin{proposition}\label{prop1}
Let $G$ be a Kaluza-Klein metric on $TM$, then the corresponding Levi-Civita connection $\bar\nabla$ is characterized, at $(p,e)\in TM$, by
\begin{align*}
\bar{\nabla}_{X^h} Y^h &= (\nabla_{X}Y)^h - \tfrac{A'}{B + |e|^2C} g(X,Y) e^v -\tfrac{1}{2} (R(X,Y)e)^v ;\\
\bar{\nabla}_{X^h} Y^v &= ( \tfrac{-B}{2A}R(Y,e)X + \tfrac{A'}{A} g(Y,e)X)^h + (\nabla_{X}Y)^v ;\\
\bar{\nabla}_{X^v} Y^h &= ( \tfrac{B}{2A}R(e,X)Y + \tfrac{A'}{A} g(X,e)Y)^h ;\\
\bar{\nabla}_{X^v} Y^v &= \tfrac{B'}{B}(g(X,e)Y^v + g(Y,e)X^v) + (C' -\tfrac{2B'C}{B})\tfrac{1}{B+|e|^2 C}g(X,e)g(Y,e)e^v \\
&+ \tfrac{C-B'}{B+ |e|^2C}g(X,Y)e^v ,
\end{align*}
all functions being evaluated at $|e|^2$ and prime denotes derivation.
\end{proposition}

\begin{remark}
Note that since $\bar{\nabla}_{X^v} Y^v$ is vertical, the fibres of $TM$ are totally geodesic.\\
The geometry of a sub-class of Kaluza-Klein metrics, called generalized Cheeger-Gromoll metrics, is studied in~\cite{BLW3}.
\end{remark}

\section{Harmonic maps}

The energy of a smooth map $\phi : (M,g) \to (N,h)$ between Riemannian manifolds is
$$ E(\phi) = \tfrac{1}{2} \int_{M} |d\phi|^2 \, v_{g} ,$$
where $|d\phi|$ is the Hilbert-Schmidt norm of $d\phi$. If $M$ is not compact, $E(\phi)$ is defined over compact subsets. Critical points of this functional are called harmonic maps, and characterized by the vanishing of the tension field~\cite{ES}:
$$\tau(\phi) = \trace \nabla d\phi =0.$$
Harmonic maps generalize not only harmonic functions but also geodesics and holomorphic maps between K{\"a}hler manifolds. The starting point of the theory is the Eells-Sampson existence result.

\begin{theorem}~\cite{ES}
Let $(M,g)$ and $(N,h)$ be compact manifolds with $\riem^{N}$ negative. Then in each homotopy class there exists a harmonic map from $(M,g)$ to $(N,h)$.
\end{theorem}

Though it has developed into a rich subject (cf.~\cite{EL1,EL2,EL3}), this theory does not lend itself to the study of vector fields. Since the tangent bundle of an $n$-dimensional manifold $M$, is itself a $2n$-dimensional manifold, one can see vector fields as maps from $M$ to $TM$, and once a Riemannian metric has been chosen on $M$, equip $TM$ with a Riemannian metric of its own. Given the canonical decomposition~\eqref{eq1} of the bitangent space and the isomorphisms between the horizontal and vertical spaces, and the tangent space of $M$, the simplest possible construction of a Riemannian metric on $TM$ is the Sasaki metric:
\begin{align*}
g_{\mathrm{Sasaki}}(X^h , Y^h) &= g(X,Y)  ;\\
g_{\mathrm{Sasaki}}(X^h , Y^v) &= 0 ;\\
g_{\mathrm{Sasaki}}(X^v , Y^v) &= g(X,Y) ,
\end{align*}
for all vectors $X,Y \in T_{p}M$, $(p,e) \in TM$ and $p\in M$. Unfortunately, elementary computations based on Proposition~\ref{prop1}, show that, for a vector field $\sigma : M \to TM$, the vertical part of $\tau(\sigma)$ is
$$\nabla^{*}\nabla \sigma =0 ,$$
and a mere integration by parts implies that $\sigma$ must be parallel, with all the topological obstructions that this implies (cf.~\cite{Nouhaud,Ishihara}).\\
An alternative metric on $TM$, was proposed by Cheeger and Gromoll (\cite{CG}), and explicated by Tricerri and Muso (\cite{MT})
\begin{align*}
g_{\mathrm{CG}}(X^h , Y^h) &= g(X,Y)  ;\\
g_{\mathrm{CG}}(X^h , Y^v) &= 0 ;\\
g_{\mathrm{CG}}(X^v , Y^v) &= \omega(e) (g(X,Y) + g(X,e)g(Y,e)) ,
\end{align*}
where $\omega(e)= \tfrac{1}{1 + |e|^2}$.\\
While this metric proved useful for other problems, it carries the same rigidity as the Sasaki metric and, when $M$ is compact, no non-parallel harmonic section or map can exist for this metric (\cite{Oniciuc,BLW1}).\\
In a first attempt to relax existence conditions, this metric was generalized in \cite{BLW1} by introducing a two-parameter family of metrics, which includes 
$g_{\mathrm{Sasaki}}$ and $g_{\mathrm{CG}}$
\begin{align*}
g_{m,r}(X^h , Y^h) &= g(X,Y)  ;\\
g_{m,r}(X^h , Y^v) &= 0 ;\\
g_{m,r}(X^v , Y^v) &= \omega^{m}(e) (g(X,Y) + rg(X,e)g(Y,e)) .
\end{align*}
Depending on the choice of $(m,r)$, $r$ positive to ensure positive definiteness, one can obtain  new harmonic maps from vector fields, e.g. the Hopf vector field from $(\sn^3 ,g_{\mathrm{can}})$ into $(T\sn^3 ,g_{2,0})$. However, in some cases, for example $\sn^2$, rigidity persists and wider classes of metrics on $TM$ are now investigated (\cite{ACP}).

Vector fields allow a richer situation than the general case of maps, since the energy functional can define several variational problems. First, harmonic maps, that is critical points of $E$ with respect to all possible variations of the map; second, harmonic sections, i.e. critical points of $E$ only with respect to variations through vector fields, and finally, topology permitting, unit harmonic sections, when variations are restricted to unit vector fields. These problems clearly sit one inside the other and when the canonical projection $\pi : TM \to M$ is a Riemannian submersion, their associated Euler-Lagrange equations can be deduced from the tension field. Indeed, the characterizing equation of harmonic sections is precisely the vertical part of the tension field (hence the same rigidity for harmonic maps and sections for the Sasaki and Cheeger-Gromoll metrics), and, for unit harmonic sections, it is the proportionality of the vertical part of the tension field and the section itself.\\
Note that C.M.~Wood extended in~\cite{Wood}, Eells-Sampson's flow technique to deduce a similar existence result for sections into a fibre bundle equipped with a Kaluza-Klein type metric (which will be our case), but the homotopy of the space of vector fields is trivial.

\begin{proposition}\label{prop2}
Let $(M,g)$ be an $m$-dimensional Riemannian manifold and $G$ a Kaluza-Klein metric on $TM$. The tension field of a vector field $\sigma : (M,g) \to (TM,G)$ is given by 
$\tau(\sigma) = [\tau^{h}(\sigma)]^h + [\tau^{v}(\sigma)]^v$, where
\begin{align*}
\tau^{h}(\sigma) &= \tfrac{-B}{A} R(\nabla_{e_{i}}\sigma,\sigma)e_{i} + \tfrac{2A'}{A} g(\nabla_{e_{i}}\sigma,\sigma)e_{i}  ;\\
\tau^{v}(\sigma) &= - \nabla^{*}\nabla \sigma + \tfrac{2B'}{B} \nabla_{X(\sigma)}\sigma \\
& + \tfrac{1}{B + |\sigma|^2 C} \left( - m A' + (C' -\tfrac{2B'C}{B})|X(\sigma)|^2 + (C-B') |\nabla \sigma|^2\right) \sigma ,
\end{align*}
where $\{e_{i}\}_{i=1,\dots,m}$ is a local orthonormal frame of $(M,g)$, $R$ its Riemann curvature tensor and $X(\sigma) = \grad \tfrac{|\sigma|^2}{2}$. 
\end{proposition}

\begin{proof}
Let $\sigma : M \to TM$ be a vector field then $d\sigma : TM \to TTM$ and, from the definition of the horizontal and vertical lifts, 
$$ d\sigma(X) = X^h + (\nabla_{X}\sigma)^v .$$
Combined with Proposition~\ref{prop1}, this yields that, in a local orthonormal frame $\{e_{i}\}_{i=1,\dots,m}$
\begin{align*}
& \tau(\sigma) = \trace \nabla d\sigma \\
&= \sum_{i=1}^{m} \bar{\nabla}_{d\sigma(e_{i})}d\sigma(e_{i}) - d\sigma(\nabla_{e_{i}}e_{i}) \\
&= \sum_{i=1}^{m} \bar{\nabla}_{e^{h}_{i} + (\nabla_{e_{i}}\sigma)^v} (e^{h}_{i} + (\nabla_{e_{i}}\sigma)^v) 
- (\nabla_{e_{i}}e_{i})^h - (\nabla_{\nabla_{e_{i}}e_{i}}\sigma)^v \\
&= \sum_{i=1}^{m} (\nabla_{e_{i}}e_{i})^h - \tfrac{A'}{B+ |\sigma|^2 C} \sigma^v + (\nabla_{e_{i}}\nabla_{e_{i}}\sigma)^v 
+ ( \tfrac{-B}{2A} R(\nabla_{e_{i}}\sigma,\sigma)e_{i} + \tfrac{A'}{A} g(\nabla_{e_{i}}\sigma,\sigma)e_{i})^h \\
&+ \tfrac{A'}{A}g(\nabla_{e_{i}}\sigma,\sigma)e_{i} + \tfrac{B}{2A} R(\sigma,\nabla_{e_{i}}\sigma)e_{i})^h + \tfrac{2B'}{B} g(\nabla_{e_{i}}\sigma,\sigma)(\nabla_{e_{i}}\sigma)^v + (C' -\tfrac{2B'C}{B})\tfrac{g^{2}(\nabla_{e_{i}}\sigma,\sigma)}{B+|\sigma|^2 C}\sigma^v\\
&+ \tfrac{C-B'}{B+ |\sigma|^2C} g(\nabla_{e_{i}}\sigma,\nabla_{e_{i}}\sigma)\sigma^v- (\nabla_{e_{i}}e_{i})^h - (\nabla_{\nabla_{e_{i}}e_{i}}\sigma)^v \\
&= \left[\sum_{i=1}^{m} \tfrac{B}{A} R(\sigma,\nabla_{e_{i}}\sigma)e_{i} + 2\tfrac{A'}{A}g(\nabla_{e_{i}}\sigma,\sigma)e_{i}\right]^h\\
&+ \left[- \nabla^{*}\nabla\sigma + \tfrac{2B'}{B} \nabla_{X(\sigma)}\sigma + \tfrac{1}{B + |\sigma|^2 C} ( -mA' + (C' - \tfrac{2B'C}{B} |X(\sigma)|^2 + (C-B')|\nabla\sigma|^2) \sigma\right]^v,
\end{align*}
hence the expression for $\tau(\sigma)$.
\end{proof}

\begin{remark}
i) Parallel vector fields are harmonic maps if and only if $A' = 0$. Therefore, from now on, with the exception of the last section, we will assume $A$ to be a constant function and harmonicity will then be independent of the value of $A$, so we can choose it to be equal to $1$.\\
ii) With this choice $A = 1$, the horizontal part satisfies $\tau^{h}(f\sigma)=f\tau^{h}(\sigma)$, for any function $f$. Moreover, harmonic sections are characterized by $\tau^{v}(\sigma)=0$ and $\tau^{h}(\sigma)$ is independent of $B$ and $C$.\\
iii) The canonical projection $\pi : (TM,G) \to (M,g)$ becomes a harmonic morphism, i.e. pulls back harmonic maps onto harmonic maps.
\end{remark}

\section{Rigidity conditions}

As the horizontal part of $\tau(\sigma)$ depends only on the geometry of $(M,g)$ and not on the choice of the functions $B$ and $C$, we can find obstructions to the existence of non-trivial harmonic maps (i.e. $\tau(\sigma)=0$).

\begin{proposition}
If $(M,g)$ is a Riemannian manifold of constant sectional curvature $\kappa$ then
$$\tau^{h}(\sigma) = \kappa ( \nabla_{\sigma}\sigma - (\Div{\sigma})\sigma) .$$
If $M$ is compact with $\kappa\neq 0$ and $\sigma$ is a harmonic map of constant length $k$ then $\Div{\sigma} =\nabla_{\sigma}\sigma =0$. Moreover
\begin{enumerate}
\item if $|\nabla\sigma|^2 \leq \kappa k^2$ then $\sigma$ is a Killing vector field;
\item if $-2\kappa k^2 \geq |L_{\sigma}g|^2$ then $\sigma$ is parallel.
\end{enumerate}
\end{proposition}

\begin{proof}
When $(M,g)$ has constant sectional curvature, the expression of $\tau^{h}(\sigma)$ is a direct consequence of Proposition~\ref{prop2}. If $\sigma$ is a harmonic map of constant norm $k$, taking the inner product of the equality $(\Div \sigma) \sigma = \nabla_{\sigma}\sigma$, gives $\sigma(\tfrac{k^2}{2}) = (\Div{\sigma})k^2$, so that $\Div{\sigma} =0$, hence $\nabla_{\sigma}\sigma =0$.\\
If $M$ is compact, we can use the Yano formula
$$\int_{M} \langle \nabla^{*}\nabla \sigma , \sigma \rangle - \ricci(\sigma,\sigma) - \tfrac{1}{2} |L_{\sigma}g|^2 + (\Div \sigma)^2 \, v_{g} =0 ,$$
to deduce that 
$$\int_{M} |\nabla \sigma|^2 - \kappa k^2 - \tfrac{1}{2} |L_{\sigma}g|^2 \, v_{g} =0 ,$$
and obtain the second part of the proposition.
\end{proof}

On surfaces, a topological obstruction appears.

\begin{proposition}\label{pp1}
Let $(M^2 , g)$ be a Riemannian surface with Gaussian curvature $K_g$. If a vector field is a harmonic map from $M^2$ to $TM^2$
 then it must vanish on the set $A = \{ x \in M^2 : K_{g}(x) \neq 0 \}$.\\
In particular, no non-zero vector field of $\sn^2$ can be a harmonic map, whatever the metric chosen on $\sn^2$.
\end{proposition}
\begin{proof}
Let $g$ be a Riemannian metric on $M^2$ and let $K_g$ be its Gaussian curvature.
First recall that if $\{X,Y\}$ is a local orthonormal frame on $M^2$ then
\begin{align*}
&\Div[\Div (X)X - \nabla_{X}X] \\
&= \langle \nabla_{X} (\langle \nabla_{Y}X , Y\rangle X - \nabla_{X}X) , X \rangle +
 \langle \nabla_{Y} (\langle \nabla_{Y}X , Y\rangle X - \nabla_{X}X) , Y \rangle \\
&= \langle \nabla_{X}\nabla_{Y} X ,Y\rangle + \langle \nabla_{Y} X , \nabla_{X}Y\rangle -  \langle \nabla_{X}\nabla_{X}X,X\rangle
 + \langle \nabla_{Y} X ,Y\rangle^2 - \langle \nabla_{Y}\nabla_{X}X,Y\rangle \\
&= - [ \langle \nabla_{Y}\nabla_{X} X ,Y\rangle - \langle \nabla_{X}\nabla_{Y} X ,Y\rangle - \langle \nabla_{\nabla_{Y} X}X ,Y\rangle]
- \langle \nabla_{X}\nabla_{X} X ,X\rangle ,
\end{align*}
but
\begin{align*}
- \langle \nabla_{X}\nabla_{X} X ,X\rangle &= |\nabla_{X}X|^2 \\
&= \langle \nabla_{X} X ,Y\rangle^2 \\
&= - \langle \nabla_{X}Y,X\rangle \langle \nabla_{X} X ,Y\rangle\\
&= - \langle \nabla_{\nabla_{X}Y}X , Y \rangle .
\end{align*}
Therefore
\begin{align}
&\Div[\Div (X)X - \nabla_{X}X] \notag \\
&= -[ \langle \nabla_{Y}\nabla_{X} X ,Y\rangle - \langle \nabla_{X}\nabla_{Y} X ,Y\rangle
- \langle \nabla_{\nabla_{Y} X}X ,Y\rangle \notag\\
&\quad - \langle \nabla_{\nabla_{X}Y}X , Y \rangle]\notag \\
&= -K_g . \label{eqp}
\end{align}
Thus, if a vector field $\sigma : (M^2,g) \to (TM^2, h_{p,q})$ is a harmonic map then necessarily the horizontal part of its tension field must vanish:
$$ K_g (\Div(\sigma)\sigma - \nabla_{\sigma}\sigma) =0.$$
Assume $\sigma$ is non-zero and let $U$ be an open subset of $A\subset M^2$ where $\sigma$ does not vanish.
First, observe that, on $U$, we have
\begin{align*}
&\Div\left(\frac{\sigma}{|\sigma|}\right)\frac{\sigma}{|\sigma|} - \nabla_{\frac{\sigma}{|\sigma|}}\frac{\sigma}{|\sigma|} \\
&=\left(\left(\frac{1}{|\sigma|}\right)\Div{\sigma} + \sigma\left(\frac{1}{|\sigma|}\right)\right) \frac{\sigma}{|\sigma|} - \frac{1}{|\sigma|^2}\nabla_{\sigma}\sigma
- \frac{1}{|\sigma|}\sigma\left(\frac{1}{|\sigma|}\right)\sigma\\
&= \frac{1}{|\sigma|^2} \left(\Div(\sigma)\sigma - \nabla_{\sigma}\sigma\right).
\end{align*}
Thus Formula~\eqref{eqp} applied to $\frac{\sigma}{|\sigma|}$ contradicts the hypothesis.
\end{proof}

When $M$ is compact the Divergence Theorem and the Kato Inequality put constraints on the vanishing of the vertical part of the tension field.
\begin{proposition}
Let $(M,g)$ be a compact Riemannian manifold and $\sigma : (M,g) \to (TM,G)$ a harmonic section , i.e. $\tau^{v}(\sigma) =0$. If $2B' + |\sigma|^2 C' \geq 0$ and $-B + |\sigma|^2 B' + |\sigma|^4 C' \leq 0$ then the norm of $\sigma$ is constant. If moreover $B + |\sigma|^2 B'$ does not vanish, then $\sigma$ must be parallel.
\end{proposition}

\begin{proof}
If $\sigma : (M,g) \to (TM,G)$ is a harmonic section then
$$ \nabla^{*}\nabla \sigma - \tfrac{2B'}{B} \nabla_{X(\sigma)}\sigma =
 \tfrac{1}{B + |\sigma|^2 C} \left( - m A' + (C' -\tfrac{2B'C}{B})|X(\sigma)|^2 + (C-B') |\nabla \sigma|^2\right) \sigma ,$$
and taking the inner-product with $\sigma$ implies
\begin{equation}\label{eq2}
\Delta \tfrac{|\sigma|^2}{2} = \tfrac{1}{B + |\sigma|^2 C} \left( -(B + |\sigma|^2 B')|\nabla\sigma|^2 + (2B'+ |\sigma|^2 C') |X(\sigma)|^2\right). 
\end{equation}
If $2B' + |\sigma|^2 C'$ is positive, then we can use the inequality $|X(\sigma)|^2 \leq |\nabla\sigma|^2|\sigma|^2$ to obtain
$$\Delta \tfrac{|\sigma|^2}{2} \leq \tfrac{1}{B + |\sigma|^2 C} \left( -B + |\sigma|^2 B' + |\sigma|^4 C' \right)|\nabla\sigma|^2,$$
and the second condition of the proposition implies that $\Delta \tfrac{|\sigma|^2}{2} \leq 0$, therefore $|\sigma|$ is constant and the left-hand side of Equation~\eqref{eq2} vanishes. 
\end{proof}

\begin{proposition}
Let $(M,g)$ be a complete Riemannian manifold with positive Ricci curvature. If $\sigma : (M,g) \to (TM,G)$ is a harmonic map and $B+|\sigma|^2B'=2B' + |\sigma|^2 C'=0$ then $\sigma$ has constant norm.
\end{proposition}
\begin{proof}
Clearly from the previous proof, $\Delta \tfrac{|\sigma|^2}{2} =0$ and, by \cite{Yau}, we conclude that $|\sigma|$ is constant.
\end{proof}

\section{Constant norm}

Vector fields of constant norm usually do not exist, but when they do, they provide particularly interesting examples and are often linked to other geometric structures. Building on the formulas of Proposition~\ref{prop2}, we can rule out their harmonicity for some combinations of the functions $B$ and $C$.

\begin{proposition}
If a vector field $\sigma : (M,g) \to (TM,G)$ of constant norm $k$, is harmonic (either as a section or a map), from a Riemannian manifold into its tangent bundle equipped with a Kaluza-Klein metric then either $\sigma$ is parallel or 
$$ B(k^2) + k^2 B'(k^2)=0.$$
\end{proposition}

\begin{proof}
If $\sigma$ has constant norm $k$, then the vertical part of its tension field becomes
\begin{equation}\label{eq3}
\nabla^{*}\nabla \sigma = \tfrac{C-B'}{B + |\sigma|^2C} |\nabla\sigma|^2 \sigma ,
\end{equation}
and taking the inner-product with $\sigma$ yields
$$|\nabla\sigma|^2 = \tfrac{C-B'}{B + |\sigma|^2C} k^2 |\nabla\sigma|^2 .$$
If $\sigma$ is not parallel, then 
$$ B(k^2) + k^2 B'(k^2)=0.$$
\end{proof}

This proposition leads to a rewriting of the condition of harmonic section.

\begin{proposition}\label{prop3}
If $\sigma : M \to TM$ is a vector field of constant norm $k$, then $\sigma : (M,g) \to (TM,G)$ is a harmonic section (i.e. $\tau^{v}(\sigma) =0$), where $G$ is a Kaluza-Klein metric for the functions $B(t) = K e^{-t/k^2}$ ($K>0$) and any choice for $C$, if and only if the unit section $\frac{\sigma}{k}$ is a unit harmonic section for the Sasaki metric on $TM$.
\end{proposition}

\begin{proof}
Recall that the Euler-Lagrange equation for unit harmonic sections with respect to the Sasaki metric on $TM$, i.e. critical points of the functional $E$ for variations through unit sections, is
$$ \nabla^{*}\nabla \alpha = |\nabla\alpha|^2 \alpha ,$$
so clearly, parallel sections are unit harmonic sections for the Sasaki metric.\\
If $\sigma$ is harmonic but not parallel, by the previous proposition, the function $B$ must satisfy
$$ B(k^2) + k^2 B'(k^2)=0,$$
which is satisfied for $B(t) = K e^{-t/k^2}$ (though this is far from the only possibility).\\ 
Moreover, for this function $B$, $\tfrac{C-B'}{B + |\sigma|^2C}= 1/k^2$ and Equation~\eqref{eq3} becomes
$$\nabla^{*}\nabla \sigma = \tfrac{1}{k^2} |\nabla\sigma|^2 \sigma ,$$
that is $\alpha = \frac{\sigma}{k}$ is a unit harmonic section.
\end{proof}

\begin{proposition}
The characteristic vector field of a Sasakian manifold $(M,g)$ is a harmonic section into $(TM,G)$, where $G$ is a Kaluza-Klein metric with $B(t) = K e^{-t}$.
\end{proposition} 

\begin{remark}
As we have already remarked, the vanishing of the horizontal part of the tension field, is independent of the choice of the Kaluza-Klein metric on $TM$. For example, the Hopf vector fields on $\sn^{2n+1}$ or some unit Killing vector fields on $\widetilde{SL_{2}(\rn)}$ or the Heisenberg three-space, turn out to be also harmonic maps (cf.~\cite{BLW2}).
\end{remark}

\section{The two-torus}

Having established in Proposition~\ref{prop3} a link between unit harmonic sections for the Sasaki metric and harmonic sections for some Kaluza-Klein metrics, we can exploit results from Wiegmink~\cite{Wiegmink} to obtain a first existence result.

\begin{theorem}
For any metric $g$ on the two-torus $\tn^2$, there exist a unitary vector field $\sigma$ and a Kaluza-Klein metric $G$ such that $\sigma : (\tn^2 ,g) \to (T\tn^2 , G)$ is a harmonic section. Moreover, $\sigma$ is a harmonic map if and only if the Gaussian curvature of $(\tn^2,g)$ vanishes.
\end{theorem}

\begin{proof}
Let $\sigma$ be a unit vector field on $\tn^2$ and $u$ a function on $\tn^2$, put $\tilde{g} = e^{-u} g$ and $\tilde{\sigma} = e^{-u}\sigma$.
Then $|\tilde{\sigma}|_{\tilde{g}}^2 =1$ and the energies of these two vector fields are (up to a constant)
$$ E(\sigma) = \tfrac{1}{2}\int_{M} B(|\sigma|_{g}^2)|\nabla\sigma|_{g}^2 \, v_{g} , \quad E(\tilde{\sigma}) = \tfrac{1}{2}\int_{M} B(|\sigma|_{\tilde{g}}^2)|\tilde{\nabla}\tilde{\sigma}|_{\tilde{g}}^2 \, v_{\tilde{g}} ,$$
since $g(\nabla_{e_{i}}\sigma,\sigma) = \tilde{g}(\tilde{\nabla}_{\tilde{e}_{i}}\tilde{\sigma},\tilde{\sigma}) =0$, where $\{e_{i}\}$ and $\{\tilde{e}_{i}\}$ are orthonormal frames with respect to $g$ and $\tilde{g}$. Note that $v_{\tilde{g}} = e^{2u} v_{g}$ and
$$|\tilde{\nabla}\sigma|_{\tilde{g}}^2 = e^{-2u} \left( |\nabla\sigma|_{g}^2 + |\grad_{g} u|_{g}^2 + 2 ((\Div{\sigma})\sigma - \nabla_{\sigma}\sigma )(u)\right) .$$
Therefore
\begin{equation}\label{eq4}
E(\tilde{\sigma}) - E(\sigma) = \tfrac{1}{2}\int_{M} B(1)\left(|\grad_{g} u|_{g}^2 + 2 (\Div{\sigma})\sigma - \nabla_{\sigma}\sigma )(u))\right) \, v_{g} ,
\end{equation}
but
\begin{align*}
\Div\left(u((\Div{\sigma})\sigma - \nabla_{\sigma}\sigma)\right) &= u\Div\left((\Div{\sigma})\sigma - \nabla_{\sigma}\sigma)\right) 
 + ((\Div{\sigma})\sigma - \nabla_{\sigma}\sigma)(u) \\
&= -u K_{g} + (\Div{\sigma})\sigma - \nabla_{\sigma}\sigma)(u) ,
\end{align*}
since 
$$\Div\left((\Div{\sigma})\sigma - \nabla_{\sigma}\sigma\right) = - K_{g},$$
as we saw in the proof of Proposition~\ref{pp1}.
Besides, by the Divergence Theorem, the second integral of \eqref{eq4} vanishes and we conclude that the quantity $E(\tilde{\sigma}) - E(\sigma)$ depends only on the function $u$, not on the section $\sigma$.\\
Now, given a Riemannian metric $g$ on $\tn^2$, choose a flat metric $\tilde{g}$ conformal to $g$. This is always possible by~\cite{KW}, and, for this flat metric, take a parallel vector field $\tilde{\sigma}$. Then $\tilde{\sigma}$ is not only a harmonic section but also an absolute minimizer of the energy functional. From our relation on the energies of $\tilde{\sigma}$ and $\sigma$, we deduce that $\sigma$ must be a unit harmonic section (since Equation~\eqref{eq4} holds only for unit sections), and we know that this implies that $\sigma$ is a harmonic section from $(\tn^2 , g)$ into $(T\tn^2 ,G)$, with the function $B(t) = Ke^{-t}$ and any function $C$.
\end{proof}

\begin{proposition}
Let $\xi$ be a Killing vector field on $(\tn^2,g)$ then $\sigma=\frac{\xi}{|\xi|}$ is a harmonic section from $(\tn^2,g)$ into $(T\tn^2,G)$ (or from an open subset of $\tn^2$ if $\xi$ vanishes at some points), where $G$ is Kaluza-Klein metric with the function $B(t) = Ke^{-t}$ and any function $C$.
\end{proposition}

\begin{proof}
Recall that the condition on $\sigma$ is $ \nabla^{*}\nabla \sigma = |\nabla\sigma|^2 \sigma.$\\
Since
\begin{align*}
\nabla\frac{1}{|\xi|} &= - \frac{\nabla|\xi|}{|\xi|^2} ; \quad \Delta\frac{1}{|\xi|} = -\frac{\Delta|\xi|}{|\xi|^2} -2 \frac{|\nabla|\xi||^2}{|\xi|^3};\\
\nabla^{*}\nabla \sigma &= K_{g} \frac{\xi}{|\xi|} + \Delta(\frac{1}{|\xi|})\xi -2 \nabla_{\nabla \tfrac{1}{|\xi|}}\xi \\
&= K_{g} \frac{\xi}{|\xi|} + ( -\frac{\Delta|\xi|}{|\xi|^2} -2 \frac{|\nabla|\xi||^2}{|\xi|^3})\xi + \frac{2}{|\xi|^2} \nabla_{\nabla |\xi|}\xi ;
\end{align*}
and
$$|\nabla \sigma|^2 = \frac{|\nabla \xi|^2}{|\xi|^2} + \frac{|\nabla|\xi||^2}{|\xi|^2 } -\frac{2}{|\xi|^3} \langle \xi,\nabla_{\nabla |\xi|}\xi\rangle ,$$
we have
\begin{align*}
\nabla^{*}\nabla \sigma - |\nabla\sigma|^2 \sigma &=K_{g} \frac{\xi}{|\xi|}  -\frac{\Delta|\xi|}{|\xi|^2}\xi -2 \frac{|\nabla|\xi||^2}{|\xi|^3}\xi + \frac{2}{|\xi|^2} \nabla_{\nabla |\xi|}\xi \\
&-\frac{|\nabla \xi|^2}{|\xi|^2}\frac{\xi}{|\xi|} - \frac{|\nabla|\xi||^2}{|\xi|^2 }\frac{\xi}{|\xi|}+ \frac{2}{|\xi|^3} \langle \xi,\nabla_{\nabla |\xi|}\xi\rangle\frac{\xi}{|\xi|}\\
&= -2 \frac{|\nabla|\xi||^2}{|\xi|^3}\xi + \frac{2}{|\xi|^2} \nabla_{\nabla |\xi|}\xi - 2\frac{|\nabla|\xi||^2}{|\xi|^2 }\frac{\xi}{|\xi|} + \frac{2}{|\xi|^3} \langle \xi,\nabla_{\nabla |\xi|}\xi\rangle\frac{\xi}{|\xi|}\\
&= \frac{2}{|\xi|^3}\left( \nabla_{\nabla \tfrac{|\xi|^2}{2}}\xi -|\nabla|\xi||^2 \xi\right) ,
\end{align*}
since 
$$K_{g}|\xi|^2=\langle \nabla^{*}\nabla \xi,\xi\rangle = |\nabla \xi|^2 + \Delta \tfrac{|\xi|^2}{2} ,$$
and
$$\frac{2}{|\xi|^3} \langle \xi,\nabla_{\nabla |\xi|}\xi\rangle\frac{\xi}{|\xi|} = \frac{2}{|\xi|^3}(\nabla |\xi|)(\tfrac{|\xi|^2}{2})\frac{\xi}{|\xi|} 
= \frac{2}{|\xi|^3}|\nabla|\xi||^2 \xi .$$
Finally, since $\xi$ is Killing, $\langle \nabla_{\xi}\xi,\xi\rangle=0$ and $\nabla \tfrac{|\xi|^2}{2} = - \nabla_{\xi}\xi$, so $\nabla_{\nabla_{\xi}\xi}\xi$ is collinear with $\xi$, that is $-\nabla_{\nabla_{\xi}\xi}\xi = f \xi$, but 
$$\langle -\nabla_{\nabla_{\xi}\xi}\xi ,\xi\rangle = |\nabla_{\xi}\xi|^2 = |\nabla\tfrac{|\xi|^2}{2}|^2 = |\xi|^2 |\nabla|\xi||^2,$$ 
hence, if $\xi$ is Killing, then $\sigma$ is harmonic.
\end{proof}

\begin{remark}
This proposition remains true on any surface.
\end{remark}

\section{Spherical vector fields}

The easiest scheme to construct vector fields on a Riemannian manifold is to take the gradient of a non-constant function. In the case of spheres, the simplest choice is a linear map given by the ambient inner-product. Unfortunately, these conformal vector fields cannot be harmonic sections (or a fortiori harmonic maps) for a Kaluza-Klein metric on $T\sn^n$.

\begin{proposition}\label{propconf}
Let $a \in \rn^{n+1}\setminus \{0\}$ ($n\geq2$) and define the function 
$$\lambda : \sn^n \to \rn ,\quad x \mapsto \lambda(x) = \langle a,x\rangle .$$ 
Then the vector field $\sigma = \grad^{\sn^n} \lambda$ is never a harmonic section (or harmonic map) for any Kaluza-Klein metric on $T\sn^n$.
\end{proposition}

\begin{proof}
If $\sigma = \grad^{\sn^n} \lambda$ with $\lambda(x) = \langle a,x\rangle$ then standard computations show that (\cite{Xin})
\begin{align*}
&\nabla^{*}\nabla \sigma = \sigma ; \quad X(\sigma) = - \lambda \sigma ; \quad \nabla_{X(\sigma)}\sigma = \lambda^2 \sigma ;\\
& |\sigma|^2= |a|^2 - \lambda^2 ; \quad |X(\sigma)|^2 = \lambda^2 (|a|^2 - \lambda^2) ; \quad |\nabla\sigma|^2 = n\lambda^2 .
\end{align*}
From these terms, one can easily check that $\tau^{v}(\sigma)=0$ if and only if
$$B - \lambda^2 (2-n)B' = \lambda^2 (|a|^2 - \lambda^2)C' + [n\lambda^2 - (|a|^2 - \lambda^2)]C .$$
But since there always exists a point $p\in \sn^n$ such that $\lambda(p) =0$, at this point the equation becomes
$$B(|a|^2) = -|a|^2 C(|a|^2) ,$$
which is impossible since this would contradict the positive definiteness of the metric $G$.
\end{proof}

While linear maps never yield harmonic maps on $\sn^n$, quadratic functions restrained to the sphere can have harmonic gradient, depending on the parity of the dimension.

\begin{proposition}
Let $\lambda$ be a quadratic form restricted to the $n$-sphere ($n\geq 2$) and $\sigma = \tfrac{1}{2} \grad^{\sn^n} \lambda$ its associated vector field. \\ If $n$ is even then $\sigma$ can never be a harmonic section from $(\sn^n , g_{\mathrm{can}})$ into $T\sn^n$ with a Kaluza-Klein metric.\\
If $n$ is odd, $\sigma$ is harmonic if and only if $\lambda$ has exactly two distinct eigenvalues $a_1 > a_2$ of multiplicity $\frac{n+1}{2}$ and $n>3$. \\
In this case, possible choices for the Kaluza-Klein metric on $T\sn^n$ are, setting $\mu= a_1 - a_2$:
\begin{itemize}
\item For $n=5$, $C=0$ and $B(t)=Ke^{\frac{-8t}{\mu^2}}$ ($K>0$) ;
\item For $n>5$, $C=0$ and (a prolongation of) 
$$B(t)=K [ \tfrac{n-3}{2} \mu^2 - (n-5)t]^{\tfrac{n+3}{n-5}}, \quad (K>0),$$
for $t\in [0,\tfrac{\mu^2}{4}]$.
\end{itemize}
\end{proposition}

\begin{proof} 
Let $\lambda$ be the restriction to $\sn^n$ of a quadratic form on $\rn^{n+1}$. Since it is symmetric, the matrix of $\lambda$ can be diagonalized into a matrix $M$ and the vector field $\sigma=\tfrac{1}{2} \nabla^{\sn^n} \lambda$ can be written 
$$\sigma(x) = M(x) - \lambda(x)x ,$$
for all $x\in \sn^n$.\\
Straightforward computations show that (cf.~\cite{BLW2})
\begin{align*}
&\nabla^{*}\nabla \sigma = (n+3)\sigma ; \quad \nabla_{X(\sigma)}\sigma = (\sigma_3 -3\lambda \sigma_2 )+ (4\lambda^2 -\lambda_2 )\sigma ;\\
& |\sigma|^2 = \lambda_2 - \lambda^2 ; \quad X(\sigma) = \sigma_2 - 2 \lambda \sigma ; \\
&|\nabla \sigma|^2 = |M|^2 -2\lambda_2 -2 \lambda \trace{M} + (n+3) \lambda^2 
\end{align*}
where for any $k\in\mathbb N^*$
$$\lambda_k(x) = M^k(x)\centerdot x, \quad \mbox{and } \sigma_k(x) = \tfrac12 \nabla\lambda_k = M^k(x)-\lambda_k(x) x,$$
and thus $\lambda_1 =\lambda$, $\sigma_1 = \sigma$.
Inspection of $\tau^{v}(\sigma)$ quickly reveals that a necessary condition for $\sigma$ to be a harmonic section is $\nabla_{X(\sigma)}\sigma$ parallel to $\sigma$, which implies that $M$ has exactly two distinct eigenvalues (cf.~\cite[Lemma 4.1]{BLW2}). Since a quadratic form with a single eigenvalue restricts to a constant function on $\sn^n$, we can assume the first eigenvalue of $M$, say $\mu$, to be strictly positive (of multiplicity $p$), and the other one to be zero.
In this case, if $x=x_\mu+x_0$ is the decomposition into eigenvectors, then (\cite{BLW2})
$$ \sigma_2 = \mu \sigma 
 ; \quad \lambda(x) = \mu |x_\mu |^2 ; \quad \lambda_2 = \mu^2 |x_\mu |^2 ,$$
and
\begin{align*}
\sigma_3 -3\lambda \sigma_2 &= (\mu^2 -3 \lambda \mu) \sigma \\
&= \mu^2 (1-3 |x_\mu|^2) \sigma .
\end{align*}
Hence
\begin{align*}
&\nabla_{X(\sigma)}\sigma = (\mu^2 -4\mu^2 |x_\mu|^2(1-|x_\mu|^2))\sigma ; \quad |\sigma|^2 = \mu^2 |x_\mu|^2(1-|x_\mu|^2) ;\\
& |X(\sigma)|^2 = \mu^2 (1-4|x_\mu|^2(1-|x_\mu|^2))|\sigma|^2 ;\\
& |\nabla \sigma|^2 = p\mu^2 -2\mu^2 (p+1)|x_\mu|^2 + (n+3)\mu^2 |x_\mu|^4 .
\end{align*}
Therefore $\tau^{v}(\sigma)=0$ if and only if
\begin{align}\label{eq5}
&(n+3)B - [2(\mu^2 - 4 |\sigma|^2) - (p -2 (p+1)|x_\mu|^2 + (n+3) |x_\mu|^4)\mu^2 ]B' \\
&= |\sigma|^2 (\mu^2 - 4|\sigma|^2)C' + [(p -2 (p+1)|x_\mu|^2 + (n+3) |x_\mu|^4)\mu^2 - (n+3)|\sigma|^2] C .\notag
\end{align}
To a point $x = x_\mu + x_0 \in \sn^n$, we associate the circle made up of points $\tilde{x} = \tilde{x}_\mu + \tilde{x}_0 \in \sn^n$ such that $|\tilde{x}_0|^2 = |x_\mu|^2$, or equivalently $|\tilde{x}_\mu|^2 = 1 - |x_\mu|^2$. Note that $|\sigma_{\tilde{x}}|^2 = |\sigma_{x}|^2$. We evaluate Equation~\eqref{eq5} at $x$ and $\tilde{x}$, and subtract to obtain
$$ (n+1-2p)(2|x_\mu|^2-1) (B' - C) = 0 , \quad \forall x\in \sn^n.$$
If $2p\neq n+1$, by continuity, this implies that $B'=C$. But then, Equation~\eqref{eq5} would become
$$(n+3)B + [-2\mu^2 + (n+11) |\sigma|^2]C = |\sigma|^2 (\mu^2 - 4|\sigma|^2)B'' ,$$
and at points where $|\sigma|^2 = \mu^2 /4$, this is
$$B(\mu^2 /4) + (\mu^2 /4)C(\mu^2 /4) =0 ,$$
which contradicts the positive definiteness of $G$.\\
Therefore, necessarily $p=\frac{n+1}{2}$, which forces $n$ to be odd and Equation~\eqref{eq5} becomes
\begin{align}\label{eq6}
&(n+3)B + [\tfrac{n-3}{2}\mu^2 - (n-5) |\sigma|^2]B' = |\sigma|^2 (\mu^2 - 4|\sigma|^2)C' + [\tfrac{n+1}{2}|\mu|^2 -2 (n+3) |\sigma|^2] C .
\end{align}
Three cases appear, $n=3$, $n=5$ and $n>5$.\\
For $n=3$, Equation~\eqref{eq6} becomes
$$6B + 2|\sigma|^2 B' = |\sigma|^2 (\mu^2 -4 |\sigma|^2)C' + 2(\mu^2 -6 |\sigma|^2)C ,$$
which can be written as an ordinary differential equation in the variable $t=|\sigma|^2 \in [0,\mu^2 /4]$:
$$3B(t) + t B'(t) = \tfrac{t}{2}(\mu^2 -4t)C'(t) + (\mu^2 -6t)C(t),$$
with $t \in [0,\mu^2 /4]$. The homogeneous equation $3B(t) + t B'(t) =0$ has $B(t) = \tfrac{K}{t^3}$ ($K\in\rn$) as solution, but it is not defined at $t=0$. If $B(t) = \tfrac{K(t)}{t^3}$ is a solution then the function $K$ must satisfy
$$K'(t) = t [\tfrac{t^2}{2} (\mu^2 -4t)C(t)]' ,$$
with $K(0)=0$. 
Therefore 
\begin{align*}
\tfrac{t^2}{2} (\mu^2 -4t)C(t) &= \int_{0}^{t} \tfrac{K'(s)}{s} \, ds \\
&= \tfrac{K(t)}{t} + \int_{0}^{t} \tfrac{K(s)}{s^2} \, ds \\
& \geq  \tfrac{K(t)}{t} = t^2 B(t),
\end{align*}
hence $(\mu^2 -4t)C(t) \geq 2 B(t)$, which, for $t=\tfrac{\mu^2}{4}$, contradicts $B$ non-negative.\\
If $n=5$, Equation~\eqref{eq6} can also be written as an ordinary differential equation in $t=|\sigma|^2$
\begin{equation}\label{eq7}
8B + \mu^2 B' = t(\mu^2 -4t)C' + (3\mu^2 -16t)C,
\end{equation}
which is satisfied for $C=0$ and $B(t) = Ke^{\frac{-8t}{\mu^2}}$ ($K>0$). For more solutions, start with $B(t) = K(t) e^{\frac{-8t}{\mu^2}}$ and then $K$ must satisfy
\begin{equation*}
K'(t)= \tfrac{1}{\mu^2} (t(\mu^2 -4t)C' + (3\mu^2 -16t)C)e^{\tfrac{8t}{\mu^2}} ,
\end{equation*}
for $t\in [0 ,\tfrac{\mu^2}{4}]$, and, for any $C$ positive (to ensure that the metric is Riemannian), we can construct a positive solution to \eqref{eq7}.\\
For $n>5$, we can take $C=0$ and the ODE becomes
$$(n+3)B = [(n-5)t - \tfrac{n-3}{2} \mu^2]B',$$
so we have the solution 
$$B(t) = K [ \tfrac{n-3}{2} \mu^2 - (n-5)t]^{\tfrac{n+3}{n-5}} , \quad (K>0),$$
for $t\in [0,\tfrac{\mu^2}{4}]$ and and we prolong $B$ to the whole of $\rn$ by continuity (keeping $B$ strictly positive).  
To obtain more solutions, put $B(t)= K(t)f(t)$, where $f(t)=[\tfrac{n-3}{2} \mu^2 - (n-5)t]^{\tfrac{n+3}{n-5}}$, then Equation~\eqref{eq6} gives the condition
$$K'(t) = [t(\mu^2 -4t)C' + (\tfrac{n+1}{2}\mu^2 -2(n+3)t)C][\tfrac{n-3}{2} \mu^2 - (n-5)t]^{\tfrac{-2(n-1)}{n-5}},$$
for $t\in [0 ,\tfrac{\mu^2}{4}]$, and we can find a primitive to construct a positive solution on $[0,\tfrac{\mu^2}{4}]$ and extend it to $\rn^{+}$.

\end{proof}

The most geometrically meaningful type of sections of the tangent bundle, is the Killing vector fields. While we can observe again different behaviour according to the parity of the dimension, harmonic Killing vector fields will exist in all cases, except dimension two.

\begin{proposition}\label{prop4}
Let $\xi$ be a Killing vector field on $\sn^{2p}$ and define its invariant axis to be
$$F_{\xi} = \ker (r_{\xi} - \mathrm{Id}),$$
where $r_{\xi}$ is the flow of $\xi$.\\
Then $\dim F_{\xi} = 2k+1$, with $0\leq k\leq p-1$ and $\xi$ is a harmonic section into $T\sn^{2p}$ equipped with a Kaluza-Klein metric $G$ if and only if
\begin{enumerate}
\item $k\neq p-1$, i.e. $\dim F_{\xi}$ is not maximal;
\item $\xi = \lambda (-x_2 , x_1, \dots ,-x_{2p-2k},x_{2p-2k-1},0,\dots ,0)$, $\lambda$ a non-zero constant.
\end{enumerate}
One possible choice for the functions defining the metric $G$ is $C=0$ and $B(t) = K e^{\tfrac{-(2p-1)t}{2\lambda^2 (p-1-k)}}$ ($K>0$).
\end{proposition}

\begin{proof}
Let $\xi$ be a Killing vector field on the even-dimensional sphere $\sn^{2p}$ and $r_{\xi}$ the isometry given by the flow of $\xi$. Its invariant axis $F_{\xi} = \ker (r_{\xi} - \mathrm{Id})$ has dimension $2k+1$ with $0\leq k\leq p-1$. Up to a permutation on the coordinates, any Killing vector field can be written $\xi = \sum_{i=1}^{p-k} \theta_{2i} \sigma_{2i}$, where the $\theta_{2i}$'s are non-zero constants and $\sigma_{2i}$ is the Killing vector field given by rotation in the $(x_{2i-1},x_{2i})$-plane. For Killing vector fields on $\sn^{2p}$, the equation of harmonic sections is
\begin{equation}\label{eq8}
(2p-1) \xi + \tfrac{2B'}{B} \nabla_{\nabla_{\xi}\xi}\xi = \tfrac{1}{B+|\xi|^2 C} \left[ (C-B')|\nabla\xi|^2 + (C' - \tfrac{2B'C}{B}) |\nabla_{\xi}\xi|^2\right] \xi,
\end{equation}
since $\nabla_{\xi}\xi= -\tfrac{1}{2} \nabla|\xi|^2$ and $\nabla^{*}\nabla \xi = (2p-1)\xi$.\\
First assume that $\theta_{2i} = \lambda$ for all $i=1,\dots,p-k$. Then $\xi = \lambda \sigma$, with $\sigma = \sum_{i=1}^{p-k}\sigma_{2i}$, and
\begin{align*}
& \nabla_{\nabla_{\xi}\xi}\xi = \lambda^3 (|\sigma|^2 -1)\sigma ; \quad |\nabla_{\xi}\xi|^2 = \lambda^4 |\sigma|^2(1-|\sigma|^2) ;\\
& |\nabla\xi|^2 = 2 \lambda^2 (p-k-|\sigma|^2),
\end{align*}
so Equation~\eqref{eq8} becomes
\begin{equation}\label{eq9}
2\lambda^2 (p-k-1)B' + (2p-1)B = \lambda^2 [2(p-k) - (2p+1)|\sigma|^2 ] C + \lambda^4 |\sigma|^2 (1-|\sigma|^2)C'.
\end{equation}
If $k\neq p-1$, then we can choose $C=0$ and $B(t) = K e^{\tfrac{-(2p-1)t}{2\lambda^2 (p-1-k)}}$ ($K>0$). For any positive function $C$, one can find a positive solution of \eqref{eq9} over the interval $[0,\lambda^2]$ and extend over to $\rn^{+}$, to obtain a suitable function $B$. \\
If $k=p-1$, Equation~\eqref{eq9} becomes
$$(2p-1)B = \lambda^4 |\sigma|^2 (1-|\sigma|^2) C' + \lambda^2 (2 - (2p+1)|\sigma|^2) C,$$
and $\xi = \lambda (-x_2 , x_1,0,\dots,0)$. Evaluated at the point $p_1 =(1,0,\dots,0)$, this equation becomes
$$(2p-1)[ B(|\xi(x)|^2) +|\xi(x)|^2 C(|\xi(x)|^2)]=0,$$
which contradicts the fact that $G$ is Riemannian.\\
Suppose now that there exist $i$ and $j$ such that $\theta_{2i} \neq \theta_{2j}$. Then
$$\xi = (-\theta_2 x_2 , \theta_2 x_1 , \dots , -\theta_{2p-2k} x_{2p-2k} , \theta_{2p-2k} x_{2p-2k-1}, 0,\dots , 0),$$
and straightforward computations show that
\begin{align*}
& \nabla_{\nabla_{\xi}\xi}\xi = \sum_{i=1}^{p-k} \theta_{2i} (|\xi|^2 -\theta_{2i}^{2} ) \sigma_{2i} ; \quad
|\nabla_{\xi}\xi|^2 = - |\xi|^4 + \sum_{i=1}^{p-k} \theta_{2i}^{4} |\sigma_{2i}|^2 ;\\
& |\nabla \xi|^2 = 2 (\sum_{i=1}^{p-k} \theta_{2i}^{2} - |\xi|^2),
\end{align*}
and Equation~\eqref{eq8} becomes
\begin{align*}
&(2p-1) \sum_{i=1}^{p-k} \theta_{2i}\sigma_{2i} - \frac{2B'}{B} \sum_{i=1}^{p-k} \theta_{2i}(\theta_{2i}^{2} - |\xi|^2) \sigma_{2i} \\
&= \frac{\sum_{i=1}^{p-k} \theta_{2i}\sigma_{2i}}{B + |\xi|^2C} \left[ (C-B') 2  (\sum_{i=1}^{p-k} \theta_{2i}^{2} - |\xi|^2) +  (C' - \tfrac{2B'C}{B}) (\sum_{i=1}^{p-k} \theta_{2i}^{4}|\sigma_{2i}|^2 - |\xi|^4)\right].
\end{align*}
From the independence of the vectors $\sigma_{2i}$, we deduce that
\begin{align*}
&(2p-1) - \tfrac{2B'}{B}( \theta_{2i}^{2} - |\xi|^2) = \tfrac{1}{B + |\xi|^2C} [ (C-B') 2  (\sum_{i=1}^{p-k} \theta_{2i}^{2} - |\xi|^2) \\
&+  (C' - \tfrac{2B'C}{B}) (\sum_{i=1}^{p-k} \theta_{2i}^{4}|\sigma_{2i}|^2 - |\xi|^4)], \quad \forall i =1,\dots,p-k, 
\end{align*}
so necessarily
$$\tfrac{2B'}{B} \theta_{2i}^{2} = \tfrac{2B'}{B} \theta_{2j}^{2} , \quad \forall i,j =1,\dots,p-k .$$ 
Since there exist $i$ and $j$ such that $\theta_{2i}\neq \theta_{2j}$, then $B'(|\xi|^2)=0$ and the equation becomes
$$(2p-1)B = [\sum_{i=1}^{p-k} \theta_{2i}^{4}|\sigma_{2i}|^2 - |\xi|^4] C' + [2\sum_{i=1}^{p-k} \theta_{2i}^{2} - (2p+1)|\xi|^2] C .$$
Let $p_1 = (1,0\dots,0)$ then $\xi(p_1 ) = \theta_2 \sigma_{2}(p_1 )$, $|\sigma_{2}(p_1 )|^2=1$ and $|\xi(p_1 )|^2 = \theta_2^2$, and the condition becomes
$$(2p-1)B(\theta_2^2) = \left(2\sum_{i=1}^{p-k} \theta_{2i}^{2} - (2p+1)\theta_2^2 \right) C(\theta_2^2),$$
which is equivalent to
$$(2p-1)\left(B(\theta_2^2) + \theta_2^2C(\theta_2^2)\right) = \left(2\sum_{i=2}^{p-k} \theta_{2i}^{2} \right) C(\theta_2^2).$$
Since $G$ is positive definite, $B(\theta_2^2) + \theta_2^2C(\theta_2^2)$ is strictly positive, hence $C(\theta_2^2)$ is also strictly positive. Therefore $2\sum_{i=1}^{p-k} \theta_{2i}^{2} - (2p+1)\theta_2^2 >0$ (recall that $B$ is a non-negative function). Similarly, testing the equation at the point $p_j = (0,\dots,0,1,0,\dots ,0)$ ($1$ at the $2j$-th position) yields
$$2\sum_{i=1}^{p-k} \theta_{2i}^{2} - (2p+1)\theta_{2j}^2 >0, \quad \forall j=1,\dots,p-k,$$
and, summing up these inequalities, we obtain
$$ -(2k+1) \sum_{i=1}^{p-k} \theta_{2i}^{2} >0 ,$$
which is clearly impossible. Hence if $\xi$ is a harmonic section, then $\theta_{2i}=\theta_{2j}$ for all $i,j=1,\dots,p-k$.
\end{proof}

\begin{remark}
Since, when $p=1$, $k$ can only take the value $0$, there exists no harmonic Killing vector field on $\sn^2$, whatever the Kaluza-Klein metric on $T\sn^2$.\\
The harmonic sections obtained in Proposition~\ref{prop4} cannot be harmonic maps because, though they are certainly divergence free, their norms cannot be constant, unless null, as Killing vector fields on even-spheres must vanish at some points.
\end{remark}

On odd-dimensional spheres, the situation is slightly different.

\begin{proposition}\label{prop5}
Let $\xi$ be a Killing vector field on $\sn^{2p+1}$ and 
$$F_{\xi} = \ker (r_{\xi} - \mathrm{Id}),$$
its invariant axis, where $r_{\xi}$ is the flow of $\xi$.\\
Then $\dim F_{\xi} = 2k$, with $0\leq k\leq p$ and $\xi$ is a harmonic section into $T\sn^{2p}$ equipped with a Kaluza-Klein metric $G$, if and only if
\begin{enumerate}
\item $k\neq p$, i.e. $\dim F_{\xi}$ is not maximal;
\item $\xi = \lambda (-x_2 , x_1, \dots ,-x_{2p-2k+2},x_{2p-2k+1},0,\dots ,0)$.
\end{enumerate}
If $k=0$, there is no condition on the function $C$ (except that the metric $G$ must be Riemannian) while $B$ must satisfy
$$\lambda^2 B'(\lambda^2) + B(\lambda^2) =0.$$
For example, $B(t) = Ke^{-t/\lambda^2}$ or $B(t) = K(1+t)^{-(1+\tfrac{1}{\lambda^2})}$ ($K>0$) are solutions. Moreover $\xi$ will be a harmonic map.\\
If $k\neq 0$, one solution is given by $C=0$ and $B(t)= K e^{-\tfrac{pt}{\lambda^2 (p-k)}}$ ($K>0$). However, as in the even dimensional case, $\xi$ will not be a harmonic map because it does not have constant norm.
\end{proposition}

\begin{proof}
Let $\xi$ be a Killing vector field on the sphere $\sn^{2p+1}$ and $r_{\xi}$ the isometry given by the flow of $\xi$. Its invariant axis $F_{\xi} = \ker (r_{\xi} - \mathrm{Id})$ has dimension $2k$ with $0\leq k\leq p$. For Killing vector fields on $\sn^{2p+1}$, the equation of harmonic sections is
\begin{equation}\label{eq10}
(2p) \xi + \tfrac{2B'}{B} \nabla_{\nabla_{\xi}\xi}\xi = \tfrac{1}{B+|\xi|^2 C} \left[ (C-B')|\nabla\xi|^2 + (C' - \tfrac{2B'C}{B}) |\nabla_{\xi}\xi|^2\right] \xi .
\end{equation}
First, assume that all the coefficients $\theta_{2i}$ are equal to $\lambda$, then $\xi = \lambda \sigma$, 
where $\sigma= \sum_{i=1}^{p-k}\sigma_{2i}$, and 
\begin{align*}
& \nabla_{\nabla_{\xi}\xi}\xi = \lambda^3 (|\sigma|^2 -1)\sigma ; \quad |\nabla_{\xi}\xi|^2 = \lambda^4 |\sigma|^2(1-|\sigma|^2) ;\\
& |\nabla\xi|^2 = 2 \lambda^2 (p+1-k-|\sigma|^2),
\end{align*}
and Equation~\eqref{eq10} becomes
\begin{equation}\label{eq11}
\lambda^2 (p-k)B' + pB = \lambda^2 [(p+1-k) - (p+1)|\sigma|^2 ] C + \frac{\lambda^4}{2}  |\sigma|^2 (1-|\sigma|^2)C'.
\end{equation}
In the special case when $\dim F_{\xi} =0$, then $\xi = \lambda \sigma$, where \\
$\sigma(x)= (-x_2 ,x_1,\dots, -x_{2p+2} , x_{2p+1})$, and $\xi$ has constant norm $|\lambda|$. Equation~\eqref{eq11} then simplifies to
$$\lambda^2 B'(\lambda^2) + B(\lambda^2)=0 ,$$
so we can choose $B(t)= K e^{-\tfrac{t}{\lambda^2}}$ ($K>0$) or $B(t) = K (1 + t)^{-(1 + \tfrac{1}{\lambda^2})}$, and $C$ to be any function (as long as the resulting metric $G$ is positive definite). Since $\xi$ has constant norm, it turns out to be also a harmonic map.\\
If $\dim F_{\xi} = 1,\dots,p-1$, then we have the solution $C=0$ and \\
$B(t)= K e^{\tfrac{-pt}{\lambda^2(p-k)}}$ ($K>0$). As in the even-dimensional case, to any positive function $C$, one can find a positive function $B$ solving \eqref{eq11}.\\
If $\dim F_{\xi} = 2p$, then the equation becomes
$$pB = \lambda^2 [1 - (p+1)|\sigma|^2 ] C + \frac{\lambda^4}{2}  |\sigma|^2 (1-|\sigma|^2)C',$$
and $\xi= \lambda (-x_2 , x_1, 0 ,\dots,0)$. We evaluate this equation at the point $p_1= (1,0,\dots,0)$ to obtain the condition
$$pB(\lambda^2) = -p\lambda^2 C(\lambda^2),$$
which contradicts the positive definiteness of the metric $G$.\\
In case there exist coefficients $i$ and $j$ such that $\theta_{2i}\neq \theta_{2j}$, $\dim F_{\xi}$ must be different from $2p$ and the decomposition
$$\xi = \sum_{i=1}^{p+1-k} \theta_{2i} \sigma_{2i},$$
leads to
\begin{align*}
&|\nabla_{\xi}\xi|^2 = -|\xi|^4 + \sum_{i=1}^{p+1-k} \theta_{2i}^{4} |\sigma_{2i}|^2 ;\quad
\nabla_{\nabla_{\xi}\xi}\xi = \sum_{i=1}^{p+1-k}\theta_{2i}(|\xi|^2-\theta_{2i}^2 )\sigma_{2i} ;\\
&  |\nabla\xi|^2 = 2(\sum_{i=1}^{p+1-k} \theta_{2i}^{2} -|\xi|^2) .
\end{align*}
Replacing in Equation~\eqref{eq10} and using the independence of the Killing vector fields $\sigma_{2i}$, yields
$$2p\theta_{2i} + \tfrac{2B'}{B} \theta_{2i}(|\xi|^2 - \theta_{2i}^{2}) = \tfrac{1}{B+|\xi|^2 C}[(C-B')|\nabla\xi|^2 + (C'-2\tfrac{B'C}{B})|\nabla_{\xi}\xi|^2]\theta_{2i},$$
therefore
$$\tfrac{2B'}{B} \theta_{2i}^2 =\tfrac{2B'}{B} \theta_{2j}^2, \quad\forall i,j=1,\dots,p+1-k,$$
hence the necessary condition $B'(|\xi|^2)=0$.
Under this condition, Equation~\eqref{eq10} becomes
$$2pB = [-|\xi|^4 + \sum_{i=1}^{p+1-k} \theta_{2i}^{4} |\sigma_{2i}|^2 ]C' + [2 \sum_{i=1}^{p+1-k}\theta_{2i}^2 -2(p+1)|\xi|^2]C,$$
evaluated at the point $p_1=(1,0,\dots,0)$ this gives
$$2pB(\theta_{2}^{2}) = [2\sum_{i=1}^{p+1-k}\theta_{2i}^2 -2(p+1)\theta_{2}^{2}]C(\theta_{2}^{2}),$$
which is equivalent to
$$2p\left(B(\theta_{2}^{2}) + \theta_{2}^{2}C(\theta_{2}^{2})\right) = [2\sum_{i=2}^{p+1-k}\theta_{2i}^2 ]C(\theta_{2}^{2}),$$
and this implies that $C(\theta_{2}^{2})$ is strictly positive and therefore, since $B$ is a strictly positive function,
$$\sum_{i=1}^{p+1-k}\theta_{2i}^2 -(p+1)\theta_{2}^2 >0 .$$
Similarly, evaluating the equation at the point $p_j = (0,\dots,0,1,0,\dots ,0)$ ($1$ at the $2j$-th position) yields
$$\sum_{i=1}^{p+1-k} \theta_{2i}^{2} - (p+1)\theta_{2j}^2 >0, \quad \forall j=1,\dots,p+1-k,$$
and summing up these inequalities we obtain
$$ -k \sum_{i=1}^{p+1-k} \theta_{2i}^{2} >0 ,$$
which is clearly impossible. Hence, $\xi$ cannot be a harmonic section.
\end{proof}

\section{A larger class of metrics}

We can render harmonic some of the vector fields considered previously, if we enlarge our class of metrics on the tangent bundle by dropping the condition $A\equiv 1$. Then the conditions for harmonic sections and harmonic maps are given by Proposition~\ref{prop2}.

\begin{proposition}
Let $a\in \rn^{n+1}\setminus\{0\}$ and $\lambda : \sn^n \to \rn , \, \lambda(x) = \langle a,x\rangle$. Then the conformal gradient field $\sigma = \grad^{\sn^n} \lambda$ is a harmonic section for 
\begin{enumerate}
\item $n=2$, $C=0$, $A(t)=B+ A_0$ ($A_0\in \rn^{+}\setminus\{0\}$) and $B(t)=K e^{-\tfrac{t}{2}}$ ($K>0$) and $\sigma$ is also a harmonic map.
\item $n>2$, $C=0$, $A(t)=B(t) + A_0$ ($A_0\in \rn^{+}\setminus\{0\}$) and $B(t)= K [ n + (n-2)|a|^2 - (n-2)t]^{\tfrac{1}{n-2}}$ ($K>0$); but $\sigma$ is not a harmonic map.
\end{enumerate}
\end{proposition}

\begin{proof}
With the notations of the proof of Proposition~\ref{propconf}, we see that a conformal gradient field $\sigma$ will be a harmonic section if
\begin{align*}
B - \lambda^2 (2-n)B' +nA'= \lambda^2 (|a|^2 - \lambda^2)C' + [n\lambda^2 - (|a|^2 - \lambda^2)]C ,
\end{align*}
and if we choose $C=0$ and $A=B+ A_0$ ($A_0 \in \rn^{+}\setminus\{0\}$), the condition becomes
$$B(t) + [(n-2)|a|^2 +n -(n-2)t] B'(t) = 0 ,$$
with $t= |\sigma|^2 = |a|^2 - \lambda^2$.\\
If $n=2$, $B(t)= K e^{-\tfrac{t}{2}}$ ($K>0$) is a solution; and if $n>2$, 
$B(t)= K [ n + (n-2)|a|^2 - (n-2)t]^{\tfrac{1}{n-2}}$ ($K>0$), is a solution, for $t\in [0,|a|^2]$, which can be extended to $\rn^{+}$.\\
For these solutions the horizontal part of the tension field for a conformal vector field is
$$-(n-1)\lambda B \sigma = 2 \lambda B' \sigma,$$
and this is only satisfied in the case $n=2$.
\end{proof}

\begin{remark}
A similar result can be found in~\cite{ACP2}.
\end{remark}

\begin{proposition}
For any Killing vector field $\xi$ on $\sn^n$ ($n\geq 2$), there exists a Kaluza-Klein metric on $T\sn^n$ making $\xi$ into a harmonic section. If $n=2$ then $\xi$ is also a harmonic map but not if $n=2p$ ($p>1$). If $n=2p+1$, $\xi$ will not be a harmonic map unless it has constant norm (cf. Proposition~\ref{prop5}).
\end{proposition}

\begin{proof}
Let $\xi : \sn^{2p} \to T\sn^{2p}$ be a Killing vector field. As in Proposition~\ref{prop4}, we call $F_{\xi} = \ker (r_{\xi} - \mathrm{Id})$ the invariant axis of its flow and $2k+1$ its dimension. Up to a permutation on the coordinates, $\sigma = \sum_{i=1}^{p-k} \theta_{2i} \sigma_{2i}$, where the $\theta_{2i}$'s are non-zero constants and $\sigma_{2i}$ is the Killing vector field given by rotation in the $(x_{2i-1},x_{2i})$-plane.
For Killing vector fields on $\sn^{2p}$, the equation of harmonic sections, with respect to the new metric $G$, is
\begin{equation}\label{eq12}
(2p-1) \xi + \tfrac{2B'}{B} \nabla_{\nabla_{\xi}\xi}\xi = \tfrac{1}{B+|\xi|^2 C} \left[ (C-B')|\nabla\xi|^2 + (C' - \tfrac{2B'C}{B}) |\nabla_{\xi}\xi|^2 -2pA' \right] \xi .
\end{equation}
First assume that the coefficients $\theta_{2i}$ are all equal to $\lambda$. If $k\neq p-1$, we know, from Proposition~\ref{prop4}, a Kaluza-Klein metric making $\xi$ a harmonic section.  If $k=p-1$ (this is the general case on $\sn^2$), Equation~\eqref{eq12} becomes
$$ (2p-1)B + 2pA' = \lambda^4 |\sigma|^2 (1-|\sigma|^2) C' + \lambda^2 (2 -(2p+1)|\sigma|^2)C,$$
where $\xi = \lambda \sigma$ (as in Proposition~\ref{prop4}).
Choosing $C=0$ and $A = B + A_0$ ($A_0 \in \rn^{+}\setminus\{0\}$), we obtain the solution $B(t) = K e^{\left(\tfrac{1}{2p}-1\right)t}$ ($K>0$).\\
In the more general case where there exist coefficients $i$ and $j$ such that $\theta_{2i}\neq \theta_{2j}$, by the same arguments as in Proposition~\ref{prop4}, we show that $B(t)=B_0$ ($B_0 \in \rn^{+}\setminus\{0\}$) and Equation~\eqref{eq12} becomes
$$(2p-1)B_0 + 2p A' = |\nabla_{\xi}\xi|^2 C' -\tfrac{1}{2}\Delta |\xi|^2 C ,$$
which admits the positive solution $A(t) = -\tfrac{2p-1}{2p} B_0 t + A_0$ for $t\in [0,\sup |\xi|^2]$ and continuously prolonged over $\rn^{+}$.\\
Similar arguments apply to the odd-dimensional case.
\end{proof}

\end{document}